\documentclass[11pt]{article}
\usepackage{amsfonts}
\usepackage{latexsym,amssymb,amsmath,graphics,cite,arydshln}
\usepackage{tikz}
\usepackage{fancyhdr}
\usepackage{lipsum}
\usepackage{lineno}
\usepackage{color}
\usepackage{caption}

\usepackage{algorithm} %format of the algorithm
\usepackage{algorithmic} %format of the algorithm
\usepackage{multirow} %multirow for format of table
\usepackage{amsmath}
\usepackage{xcolor}

\begin{document}
\newcommand{\qed}{\hphantom{.} \hfill $\Box$\medbreak}
\newcommand{\proof}{\noindent{\bf Proof \ }}

\newtheorem{theorem}{Theorem}[section]
\newtheorem{lemma}[theorem]{Lemma}
\newtheorem{corollary}[theorem]{Corollary}
\newtheorem{remark}[theorem]{Remark}
\newtheorem{example}[theorem]{Example}
\newtheorem{definition}[theorem]{Definition}
\newtheorem{fact}[theorem]{Fact}
\newtheorem{claim}[theorem]{Claim}
\newtheorem{proposition}[theorem]{Proposition}
\newtheorem{note}[theorem]{Note}

\begin{center}
{\Large\bf The area of reduced spherical polygons \footnote{Supported by NSFC under Grant $11971053$  and Science Foundation
of Hebei Normal University L2020Z01.}
}

\vskip12pt
Cen Liu $^{\dag}$, Yanxun Chang $^{\dag}$, Zhanjun Su $^{\ddag}$ \\[2ex] {$^{\dag}$\footnotesize Department of Mathematics, Beijing Jiaotong University, Beijing 100044, China

$^{\ddag}$School of Mathematical Sciences, Hebei Normal University, Shijiazhuang, 050024, China.}

{\footnotesize
18118005@bjtu.edu.cn, yxchang@bjtu.edu.cn, suzj888@163.com}
\vskip12pt
\end{center}

\vskip12pt

\noindent {\bf Abstract:}
We confirm two conjectures of Lassak on the area of reduced spherical polygons.
The area of every reduced spherical non-regular $n$-gon is less than that of the regular spherical $n$-gon of the same thickness.
Moreover, the area of every reduced spherical polygon is less than that of the regular spherical odd-gons of the same thickness and whose number of vertices tends to infinity.
\vskip12pt

\noindent {\bf Keywords}: reduced convex body, spherical polygon, thickness, area

\noindent {\bf Mathematics Subject Classification (2010)}: 52A55
%%%%%%%%%%%%%%%%%%%%%%%%%%%%%%%%%%%%%%%%%%%%%%%%%%%%%%%%%%%%%%%%%%%%%%%%%

\section{Introduction}
We focus on the reduced convex bodies introduced by Heil in~\cite{Heil}.
Reduced convex bodies are helpful for solving various extremal problems concerning the minimal width of convex bodies. Some basic properties of the reduced convex bodies in two-dimensional Euclidean space $E^{2}$ are introduced by Lassak in \cite{1990}. Lassak \cite{2005} demonstrates that in $E^{2}$ the area of every reduced non-regular $n$-gon is less than that of the regular $n$-gon of the same thickness.

The notions about reduced convex bodies are extended to the $d$-dimensional unit sphere $S^{d}$ in \cite{convexbody}.
Lassak \cite{convexbody, 2020} discusses the properties of reduced convex bodies on $S^{d}$;
he \cite{polygon} further characterizes reduced convex polygons on $S^{2}$ and proposes the following conjectures:

{(1) \it The area of every reduced spherical polygon is less than that of the regular spherical odd-gons of the same thickness and whose number of vertices tends to infinity}.

 {(2) \it The area of every reduced spherical non-regular $n$-gon is less than that of the regular spherical $n$-gon of the same thickness}.

In Section $2$, we present the necessary notions of reduced spherical convex bodies and review some results in the literature.
Several useful lemmas are established in Section $3$.
Then Section $4$ aims to confirm the above two conjectures.

\section{Preliminaries}
Let $S^{2}$ be the unit sphere in $E^{3}$ centered at the origin. In this paper,
all the notions are discussed in $S^{2}$.
 A {\it great circle} is the intersection of $S^{2}$ with any two-dimensional subspace of $E^{3}$.
 A pair of {\it antipodes} are the intersection of $S^{2}$ with any one-dimensional subspace of $E^{3}$.

Clearly, if different points $a,b\in S^{2}$ are not antipodes, then there is exactly one great circle containing them;
denote by {\it arc ab}, shortly $ab$, the shorter part of the great circle containing them.
The {\it spherical distance} $|ab|$, or shortly {\it distance} of $a$ and $b$ is the length of $ab$.

A set $C\subset S^{2}$ is called {\it convex},
if it does not contain any pair of antipodes of $S^{2}$ and for arbitrary points $a,b\in C$, it is true that arc $ab\subseteq C$.
The {\it convex body} is a closed convex set with non-empty interior.

By a {\it spherical disk} of radius $r\in(0,\pi/2]$ and {\it center} $k\in S^{2}$ we mean the set $B=\{p:|pk|\leq r, p\in S^{2}\}$;
and the boundary $\mathrm{bd}(B)$ is called a {\it spherical circle}.
The spherical disk of radius $\pi/2$ is called a {\it hemisphere}.
If hemispheres $G$ and $H$ are different and their centers are not antipodes,
then $L=G\cap H$ is called a {\it lune} of $S^{2}$.
The parts of $\mathrm{bd}(G)$ and $\mathrm{bd}(H)$ contained in $G\cap H$ are
denoted by $G/H$ and $H/G$, respectively.
 We define the {\it thickness} $\Delta(L)$ of the lune $L=G\cap H $ as the distance of the centers of $G/H$ and $H/G$.

We recall some notions in \cite{convexbody}.
We say that a hemisphere $H$ {\it supports} a convex body $C\subseteq S^{2}$ at point $p$ if $C\subseteq H$ and $p\in\mathrm{bd}(H)\cap C$.
For any convex body $C\subseteq S^{2}$ and any hemisphere $K$ supporting $C$, we define the the $width$ of $C$ {\it determined} by $K$ as the minimum thickness of a lune $K\cap K^{*}$ over all hemispheres $K^{*}\neq K$ supporting $C$ and we denote it by $\mathrm{width}_{K}(C)$;
the {\it thickness} of $C$ is defined by
$$\Delta(C)=\min \{\mathrm{width}_{K}(C): \text{K is a supporting hemisphere of C}\}.$$
The thickness of $C$ is nothing else but the minimum thickness of a lune containing $C$.
A convex body $C\subset S^{2}$ is said to be {\it reduced}
if $\Delta(R)<\Delta(C)$ for each convex body
$R$ being a proper subset of $C$.

We recall some definitions given in \cite{polygon}.
Let $p$ be a point in a hemisphere different from its center and let $l$ be the great circle bounding this hemisphere. The {\it projection} of $p$ on $l$ is the point $t$ such that $|pt|=\min \{|pc|:c\in l\}$.
If $C$ is a subset of a convex set of $S^{2}$, then the intersection of all convex sets containing $C$ is called a {\it convex hull} of $C$.
The convex hull of $k\geq 3$ points on $S^{2}$
such that each of them does not belong to the convex hull of the remaining
points is called a {\it spherically convex} $k$-gon.
If $V$ is a spherically convex $k$-gon, we denote by $v_{1},\ldots ,v_{k}$ the vertices of $V$ in the counterclockwise order.
A spherically convex polygon with sides of equal length and interior angles of equal measure is called a {\it regular spherical polygon}.

\begin{lemma}{\rm\cite[Theorem 41.2]{IP}}\label{Gir-polygon}
Let $V$ be a convex $n$-gon in the unit sphere with angles $\gamma_{1},\gamma_{2},\ldots, \gamma_{n}$,
then $\mathrm{area}(V)=\gamma_{1}+\ldots+\gamma_{n}-(n-2)\pi$.
\end{lemma}

For a convex odd-gon $V=v_{1}v_{2}\cdots v_{n}$, by the {\it opposite side to the vertex $v_{i}$} we mean the side $v_{i+(n-1)/2}v_{i+(n+1)/2}$, the indices are taken modulo $n$.

\begin{lemma}{\rm\cite[Theorem 3.2]{polygon}}\label{3.2}
Every reduced spherical polygon is an odd-gon of thickness at most $\frac{\pi}{2}$.
A spherically convex odd-gon $V$ with $\Delta(V)<\frac{\pi}{2}$ is reduced if and only if the projection of every its vertices on the great circle containing the opposite side belongs to the relative interior of this side and the distance of this vertex from this side is $\Delta(V)$.
\end{lemma}

\begin{lemma}{\rm\cite[Corollary 3.3]{polygon}}\label{3.3}
Every regular spherical odd-gon of thickness at most $\frac{\pi}{2}$ is reduced.
\end{lemma}

\section{Working lemmas}
We recall a few formulas of spherical geometry in \cite{DA} which are the basic method for the research.
Consider the right spherical triangle with hypotenuse $c$ and legs $a,b$,
we use $A,B$ and $C$ to represent the corresponding angles of edges $a,b$ and $c$, respectively.
Then
\begin{equation}\label{eq2}
\cos B=\cos b\sin A,
\end{equation}
\begin{equation}\label{eq1}
\cos A=\tan b\cot c,
\end{equation}
\begin{equation}\label{eq3}
\sin b=\sin c\sin B.
\end{equation}

\begin{figure}[htb]
\centering
  \includegraphics[width=8 cm]{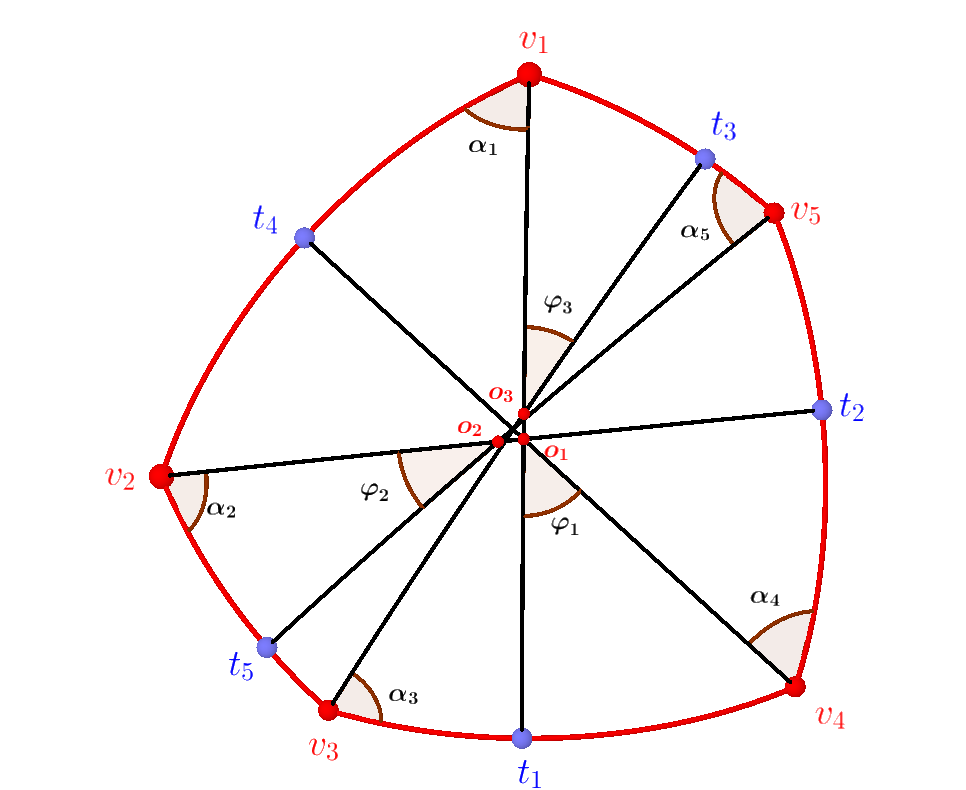}\\
\caption{Some notations}
  \label{example}
\end{figure}

In a reduced spherical polygon $V=v_{1}v_{2}\cdots v_{n}$, according to Lemma \ref{3.2}, we give some related notations.
Denote by $t_{i}$ the projection of $v_{i}$ in the opposite side $v_{i+(n-1)/2}v_{i+(n+1)/2}$.
Denote by $o_{i}$ the intersection of $v_{i}t_{i}$ and $v_{i+(n+1)/2}t_{i+(n+1)/2}$;
put $\alpha_{i}=\angle v_{i+1}v_{i}t_{i}$, $\beta_{i}=\angle t_{i}v_{i}v_{i+(n+1)/2}$,
 and $\varphi_{i}=\angle v_{i}o_{i}t_{i+(n+1)/2}=\angle t_{i}o_{i}v_{i+(n+1)/2}$, where $i\in\{1,2,\ldots,n\}$.
For example, Fig.\,\ref{example} presents some notations in a reduced spherical pentagon.

\begin{lemma}{\rm\cite[Corollary 3.9]{polygon}}\label{3.9}
If $V=v_{1}v_{2}\cdots v_{n}$ is a reduced spherical polygon with $\Delta(V)<\frac{\pi}{2}$, then $\beta_{i}\leq\alpha_{i}$.
\end{lemma}

Actually, we can gain the following lemma by Corollary $3.6$ in \cite{polygon}.
Here we prove it in a different way.

\begin{lemma}\label{congruent triangle}
For every reduced spherical polygon $V=v_{1}v_{2}\cdots v_{n}$ with $\Delta(V)<\frac{\pi}{2}$, the spherical triangles $v_{i}o_{i}t_{i+(n+1)/2}$ and $v_{i+(n+1)/2}o_{i}t_{i}$ are congruent, where $i\in\{1,2,\ldots,n\}$.
\end{lemma}
\proof
Lemma \ref{3.2} shows that
$|v_{i}t_{i}|=|v_{i+(n+1)/2}t_{i+(n+1)/2}|=\Delta(V)$, where $i\in\{1,\ldots,n\}$.
From this and $|v_{i}v_{i+(n+1)/2}|=|v_{i+(n+1)/2}v_{i}|$, we find that the right spherical triangles $v_{i}t_{i}v_{i+(n+1)/2}$ and $v_{i+(n+1)/2}t_{i+(n+1)/2}v_{i}$ are congruent.
Then we acquire
\begin{equation}\label{angle-sum}
\alpha_{i}+\beta_{i}
=\angle t_{i+(n+1)/2}v_{i}v_{i+(n+1)/2}=\angle v_{i}v_{i+(n+1)/2}t_{i}
\end{equation}
and
\begin{equation}\label{angle}
\beta_{i}=\angle t_{i}v_{i}v_{i+(n+1)/2}=\angle v_{i}v_{i+(n+1)/2}t_{i+(n+1)/2}.
\end{equation}
From \eqref{angle-sum} and \eqref{angle}, we obtain that
$\alpha_{i}=\angle t_{i+(n+1)/2}v_{i}o_{i}=\angle t_{i}v_{i+(n+1)/2}o_{i}$.
Applying this and $\angle v_{i}o_{i}t_{i+(n+1)/2}=\angle t_{i}o_{i}v_{i+(n+1)/2}=\varphi_{i}$, it turns out that the right spherical triangles
$v_{i}o_{i}t_{i+(n+1)/2}$ and $v_{i+(n+1)/2}o_{i}t_{i}$ are congruent.\hfill\qed

\begin{lemma}\label{regular}
All the vertices of a regular spherical polygon are contained in a spherical circle.
\end{lemma}
\proof
Let $V=v_{1}v_{2}\cdots v_{n}$ be a regular spherical polygon. Denote by $\gamma$ the length of the side of $V$ and by $\theta$ the interior angle of $V$.
Let $D$ be a spherical circle passing through $v_{1}$, $v_{2}$ and $v_{3}$,
and denote by $o$ be the center of the spherical disk whose boundary is $D$.
Connect $o$ with $v_{i}$, where $i\in\{1,2\ldots,n\}$.

In the spherical triangles $ov_{1}v_{2}$ and $ov_{2}v_{3}$, we have $|ov_{1}|=|ov_{2}|=|ov_{3}|$. By the definition of regular spherical polygon, we have
$|v_{1}v_{2}|=|v_{2}v_{3}|=\gamma$, then $ov_{1}v_{2}$ and $ov_{2}v_{3}$ are congruent.
Since $|ov_{1}|=|ov_{2}|=|ov_{3}|$, we have $\angle ov_{i}v_{i+1}=\angle ov_{i+1}v_{i}=\frac{\theta}{2}$, where $i=1,2$.

In the spherical triangles $ov_{1}v_{2}$ and $ov_{3}v_{4}$.
From $|ov_{1}|=|ov_{3}|$, $\angle ov_{1}v_{2}=\angle ov_{3}v_{4}=\frac{\theta}{2}$, and $|v_{1}v_{2}|=|v_{3}v_{4}|$, we obtain that $ov_{1}v_{2}$ and $ov_{3}v_{4}$ are congruent. Hence we have $|ov_{1}|=|ov_{2}|=|ov_{3}|=|ov_{4}|$.

Similarly, the spherical triangles $ov_{1}v_{2}$ and $ov_{j}v_{j+1}$ are congruent, where $j\in\{4,\ldots,n\}$ and the indices are taken modulo $n$. Hence $|ov_{1}|=|ov_{2}|=|ov_{j}|=|ov_{j+1}|$ with $j\in\{4,\ldots,n\}$.
Consequently, we have $|ov_{1}|=|ov_{2}|=\cdots=|ov_{n}|$,
and thus $v_{1},v_{2},\cdots,v_{n}$ are all contained in $D$.\hfill\qed

\begin{fact}\label{fact-varphi}
For every reduced spherical polygon $V=v_{1}v_{2}\cdots v_{n}$ with $\Delta(V)<\frac{\pi}{2}$, we have $0<\varphi_{i}<\frac{\pi}{2}$, where $i\in\{1,\ldots,n\}$.
\end{fact}
\proof
Consider the right spherical triangle $v_{i}t_{i}v_{i+(n+1)/2}$ for every $i\in\{1,\ldots,n\}$.
From \eqref{eq2}, we obtain $\cos(\alpha_{i}+\beta_{i})=\cos|v_{i}t_{i}|\sin\beta_{i}$.
 As $|v_{i}t_{i}|=\Delta(V)<\frac{\pi}{2}$, we get $\alpha_{i}+\beta_{i}<\frac{\pi}{2}$.
Applying this and Lemma \ref{3.9}, it follows that $2\beta_{i}\leq\alpha_{i}+\beta_{i}<\frac{\pi}{2}$.

It is obvious that the area of the spherical triangle $o_{i}v_{i}v_{i+(n+1)/2}$ is nonnegative. Hence by Lemma \ref{3.3} and \eqref{angle}, we get that $\beta_{i}+\beta_{i}+(\pi-\varphi_{i})-\pi\geq0$.
Therefore, $\varphi_{i}\leq2\beta_{i}<\frac{\pi}{2}$. The proof is complete.
\hfill\qed

By a {\it rotation} of a set $C\subseteq S^{2}$ around a point $p\in S^{2}$, we mean the rotation of $C$ around the straight line through $p$ and the center of $E^{3}$.
For any two points $a,b\in \mathrm{bd}(C)$, $ab$ is called a {\it chord} of $C\subseteq S^{2}$.

\begin{lemma}\label{sum}
For any reduced spherical polygon $V=v_{1}v_{2}\cdots v_{n}$ with $\Delta(V)<\frac{\pi}{2}$, we have $\sum_{i=1}^n\varphi_{i}\geq\pi$.
\end{lemma}
\proof
Lemma \ref{3.2} shows that $V$ is an odd-gon.
We use a similar technic as that of Lemma in \cite{2005} to show the statement.
Let $B_{i}=v_{i}o_{i}t_{i+(n+1)/2} \cup v_{i+(n+1)/2}o_{i}t_{i}$,
then $B_{i}\subseteq V$, where $i\in\{1,2\ldots,n\}$. Thus we have $B_{1}\cup\cdots\cup B_{n}\subseteq V$.

We intend to show $V\subseteq B_{1}\cup\cdots\cup B_{n}$.
We present every $B_{i}$ as the union of chords of $V$ which pass through $o_{i}$.
All the chords of successively $B_{1},B_{1+(n+1)/2},\ldots,B_{(n+1)/2}$ are in big circles which step by step rotate changing the centers of rotation; those centers successively are $o_{1},o_{1+(n+1)/2},\ldots,o_{(n+1)/2}$.
We assume that all the above chords in $B_{i}$ are oriented with the origins in $v_{i}t_{i+(n+1)/2}$.
For any point $p\in V$, we assume that $p$ is in the left hand side of $v_{1}t_{1}$.
When we start from $v_{1}t_{1}$, after total rotation by $\varphi_{1}+\cdots+\varphi_{n}$, we arrive at $t_{1}v_{1}$ which has the opposite direction.
Now $p$ is in the right hand side of the oriented chord $v_{1}t_{1}$.
Since the described changes of $v_{1}t_{1}$ are continuous, there is a position such that the chord contains $p$. Hence $p\in B_{1}\cup\cdots\cup B_{n}$ and then $V\subseteq B_{1}\cup\cdots\cup B_{n}$

Consequently, $V=B_{1}\cup\cdots\cup B_{n}$.

{\bf Claim 1:} If $V$ is a non-regular reduced spherical polygon, then $\sum_{i=1}^n\varphi_{i}\geq\pi$.

From Lemmas \ref{Gir-polygon} and \ref{congruent triangle}, the area of $V$ is $S_{V}=2(\alpha_{1}+\cdots+\alpha_{n})-(n-2)\pi$.
The area of $B_{i}$ is $S_{B_{i}}=2(\varphi_{i}+\alpha_{i}-\frac{\pi}{2})$.
From $V=B_{1}\cup\cdots\cup B_{n}$,
we have $S_{V}\leq\sum_{i=1}^nS_{B_{i}}$, that is
\begin{equation}\label{Sum}
2\sum_{i=1}^n\alpha_{i}-(n-2)\pi\leq2\sum_{i=1}^n(\varphi_{i}+\alpha_{i}-\frac{\pi}{2}),
\end{equation}
and then $\sum_{i=1}^n\varphi_{i}\geq\pi$.

{\bf Claim 2:} If $V$ is a regular spherical polygon, then $\sum_{i=1}^n\varphi_{i}=\pi$ and $\varphi_{i}=\frac{\pi}{n}$, where $i\in\{1,2,\ldots,n\}$.

In this case, by Lemma \ref{regular}, we have $o_{1}=\cdots=o_{n}=o$ (the notation $o$ is described in Lemma \ref{regular}).
Then, $B_{i}=v_{i}ot_{i+(n+1)/2} \cup v_{i+(n+1)/2}ot_{i}$, where $i\in\{1,2,\ldots,n\}$. Clearly, for arbitrary $i,j\in\{1,2,\ldots,n\}$, $B_{i}$ and $B_{j}$ are congruent, and their interiors satisfy $\mathrm{int}(B_{i})\cap\mathrm{int}(B_{j})=\emptyset$.
Thus we have $S_{V}=\sum_{i=1}^nS_{B_{i}}$ and then \eqref{Sum} becomes
\begin{equation*}
2\sum_{i=1}^n\alpha_{i}-(n-2)\pi=2\sum_{i=1}^n(\varphi_{i}+\alpha_{i}-\frac{\pi}{2}).
\end{equation*}
Hence we obtain that $\sum_{i=1}^n\varphi_{i}=\pi$ and $\varphi_{i}=\frac{\pi}{n}$, where $i\in\{1,2,\ldots,n\}$.\hfill\qed

In the following lemmas, we investigate the monotonicity and concavity of two kinds of functions, respectively, which are needed in Section $4$.

\begin{lemma}\label{decreasing function}
Let $f_{1}(x)=\arccos\frac{x\sqrt{1+\lambda^{2}}}{\lambda-x}$ and $f_{2}(x)=\arccos\frac{x(1+\lambda x)}{\lambda-x}$. Then $\frac{f_{1}(x)}{f_{2}(x)}$ is a decreasing function of $x$, where $\lambda\in(0,+\infty)$ and $x\in(0,\frac{-1+\sqrt{1+\lambda^{2}}}{\lambda})$.
\end{lemma}
\proof
Set $f(x)=\frac{f_{1}(x)}{f_{2}(x)}$.
Let us show that $f'(x)<0$. The derivative of $f(x)$ is
\begin{equation*}
f'(x)=\frac{\lambda\cdot h(x)}{(\lambda-x)\sqrt{1+x^{2}}\sqrt{-\lambda^{2} x^{2}-2\lambda x+\lambda^{2}}f_{2}^{2}(x)},
\end{equation*}
where
\begin{equation*}
 h(x)=-\sqrt{(1+\lambda^{2})(1+x^{2})}f_{2}(x)+(-x^{2}+2\lambda x+1)f_{1}(x).
\end{equation*}
From $x<\frac{-1+\sqrt{1+\lambda^{2}}}{\lambda}<\lambda$, and if $h(x)<0$, then $f'(x)<0$.
The first derivative of $h(x)$ is
$$h'(x)=-\frac{x\sqrt{1+\lambda^{2}}}{\sqrt{1+x^{2}}}f_{2}(x)
+2(\lambda-x)f_{1}(x).$$
The second derivative of $h(x)$ is
$$
h''(x)=-\frac{\lambda\sqrt{1+\lambda^{2}}}{\sqrt{-\lambda^{2}x^{2}-2\lambda x+\lambda^{2}}}\cdot\frac{2\lambda-3x-x^{3}}{(1+x^{2})(\lambda-x)}-\frac{\sqrt{1+\lambda^{2}}}{(1+x^{2})^{3/2}}f_{2}(x)
-2f_{1}(x).
$$
Since $x\in(0,\frac{-1+\sqrt{1+\lambda^{2}}}{\lambda})$, it follows that $\lambda-x>0$ and $2\lambda-3x-x^{3}>0$. Also we have $f_{1}(x)>0$ and $f_{2}(x)>0$.
Hence $h''(x)<0$ and then $h'(x)>h'(\frac{-1+\sqrt{1+\lambda^{2}}}{\lambda})=0$.
Thus from $h'(x)>0$, we get $h(x)<h(\frac{-1+\sqrt{1+\lambda^{2}}}{\lambda})=0$.
Therefore, $f'(x)<0$.
\hfill\qed

\begin{lemma}\label{F(z)}
 Let $F(x)=\arcsin\frac{g(x)\sqrt{1+\lambda^{2}}}{\lambda-g(x)}$
and $g(x)=\frac{-(1+\cos x)+\sqrt{(1+\cos x)^{2}+4\lambda^{2}\cos x}}{2\lambda}$,
where $\lambda\in(0,+\infty)$ and $x\in(0,\frac{\pi}{2})$. Then $F'(x)<0$ and $F''(x)<0$.
\end{lemma}
\proof
For convenience, set $r(x)=\sqrt{(1+\cos x)^{2}+4\lambda^{2}\cos x}$.
We find the first derivative of $F(x)$ is
\begin{equation*}
F'(x)=-\frac{\lambda\sqrt{2+2\lambda^{2}}\sin x}
{r(x)\sqrt{1-\cos x}\sqrt{1+2\lambda^{2}+\cos x-r(x)}}
.
\end{equation*}
Therefore, we obtain $F'(x)<0$.  The second derivative of $F(x)$ is
\begin{equation*}
F''(x)=\frac{\lambda\sqrt{2\lambda+2\lambda^{2}}\sin^{4}\frac{x}{2}
\big(-2(1+\cos x)^{2}-8\lambda^{2}+2(1+\cos x)r(x)\big)}
{(1-\cos x)^{\frac{3}{2}}r^{3}(x)\sqrt{1+2\lambda^{2}+\cos x-r(x)}}.
\end{equation*}
We can check that $-2(1+\cos x)^{2}-8\lambda^{2}+2(1+\cos x)r(x)<0$ and
$(1-\cos x)^{\frac{3}{2}}r^{3}(x)>0$ in the domain $x\in(0,\frac{\pi}{2})$. Hence $F''(x)<0$ and then $F(x)$ is a concave function of $x$.\hfill\qed

\section{The area of reduced spherical polygons }
This section aims to prove the conjectures mentioned in the introduction.
For ease of notations, we use $\omega$ to replace the thickness of a reduced spherical polygon in this part. Here we have $\omega\in (0,\pi/2)$. Let $\lambda=\tan\omega$, then $\lambda\in(0,+\infty)$.
Denote by $S$ the {\it area} of a reduced spherical polygon.

Let us define several functions which are needed in the following theorems.
Set $$f(x)=\arcsin\frac{x\sqrt{1+\lambda^{2}}}{\lambda-x},f_{2}(x)=\arccos\frac{x(1+\lambda x)}{\lambda-x},$$
$$f_{1}(x)=\frac{\pi}{2}-f(x)=\arccos\frac{x\sqrt{1+\lambda^{2}}}{\lambda-x},$$
where $x\in\big(0,(-1+\sqrt{1+\lambda^{2}})/\lambda\big)$.
Set $$g(\varphi)=\frac{-(1+\cos\varphi)+\sqrt{(1+\cos\varphi)^{2}+4\lambda^{2}\cos \varphi}}{2\lambda},$$ where $\varphi\in(0,\pi/2)$. And thus $g(\varphi)\in\big(0,(-1+\sqrt{1+\lambda^{2}})/\lambda\big)$.

\begin{lemma}\label{area of V}
For a reduced spherical polygon $V=v_{1}\cdots v_{n}$ with $\omega\in (0,\frac{\pi}{2})$, the area is $S=2\sum_{i=1}^n f(y_{i})-(n-2)\pi$,
where $y_{i}=g(\varphi_{i})$.
\end{lemma}
\proof
For every $i\in\{1,2,\dots,n\}$, we focus on the right spherical triangle $o_{i}t_{i}v_{i+(n+1)/2}$.
Put $|o_{i}t_{i}|=b_{i}$ and $|o_{i}v_{i+(n+1)/2}|=c_{i}$.
By Lemma \ref{congruent triangle} we have $|o_{i}t_{i}|+|o_{i}v_{i+(n+1)/2}|=b_{i}+c_{i}=\omega$.
Here we have $b_{i}<\omega$.
From \eqref{eq1}, we obtain
\begin{equation}\label{in1}
\cos\varphi_{i}=\frac{\tan b_{i}}{\tan c_{i}}=\frac{\tan b_{i}}{\tan(\omega-b_{i})}
=\frac{\tan b_{i}(1+\tan\omega\tan b_{i})}{\tan\omega-\tan b_{i}}.
\end{equation}
 By a simple calculation, we get that $\tan b_{i}=g(\varphi_{i})$, for simplicity, we denote $g(\varphi_{i})$ by $y_{i}$. Hence $b_{i}=\arctan y_{i}$.
Fact \ref{fact-varphi} shows that $0<\varphi_{i}<\frac{\pi}{2}$. From this and $b_{i}<\omega$, it follows that $y_{i}\in(0,\frac{-1+\sqrt{1+\lambda^{2}}}{\lambda})$.

Hence \eqref{in1} becomes $\cos\varphi_{i}=\frac{y_{i}(1+\lambda y_{i})}{\lambda-y_{i}}$
and then $\varphi_{i}=f_{2}(y_{i})=\arccos\frac{y_{i}(1+\lambda y_{i})}{\lambda-y_{i}}$.
Moreover, we have $\tan c_{i}=\frac{\tan b_{i}}{\cos\varphi_{i}}=\frac{\lambda-y_{i}}{1+\lambda y_{i}}$
and thus $c_{i}=\arctan\frac{\lambda-y_{i}}{1+\lambda y_{i}}$.

From \eqref{eq3}, we obtain
$$\sin\alpha_{i}=\frac{\sin b_{i}}{\sin c_{i}}
=\frac{\sin \arctan y_{i}}{\sin \arctan\frac{\lambda-y_{i}}{1+\lambda y_{i}}}
=\frac{y_{i}\sqrt{1+\lambda^{2}}}{\lambda-y_{i}},$$
and thus $\alpha_{i}=f(y_{i})=\arcsin\frac{y_{i}\sqrt{1+\lambda^{2}}}{\lambda-y_{i}}$.
Then Lemmas \ref{Gir-polygon} and \ref{congruent triangle} imply that the area of $V$ is
\begin{equation}\label{area}
S=2\sum_{i=1}^n \alpha_{i}-(n-2)\pi
=2\sum_{i=1}^n f(y_{i})-(n-2)\pi.
\end{equation}
\hfill\qed

\begin{theorem}\label{thm3.1}
The regular spherical $n$-gon has the maximum area among all regular spherical $k$-gons of fixed thickness, with odd numbers $k,n$ and $3\leq k\leq n$.
\end{theorem}
\proof
Let $V=v_{1}v_{2}\cdots v_{k}$ be a regular spherical odd-gon.
Lemma \ref{3.3} shows that $V$ is reduced, then  we use the same notations as that in Lemma \ref{area of V}.
By Claim $2$ in Lemma \ref{sum},
we have $\varphi_{1}=\cdots=\varphi_{k}=\frac{\pi}{k}$.
Thus $y_{1}=\cdots=y_{k}=g(\frac{\pi}{k})$, where
$$
g(\frac{\pi}{k})=\frac{-(1+\cos\frac{\pi}{k})+\sqrt{(1+\cos\frac{\pi}{k})^{2}+4\lambda^{2}\cos\frac{\pi}{k}}}{2\lambda},
$$
and $g(\frac{\pi}{k})\in(0,\frac{-1+\sqrt{1+\lambda^{2}}}{\lambda})$. For simplicity, we denote $\varphi_{i}$ and $y_{i}$ by $\varphi$ and $y$, respectively, for every $i\in\{1,2,\ldots,k\}$.
From the proof process of Lemma \ref{area of V}, we can easily see that $k=\frac{\pi}{\varphi}=\frac{\pi}{f_{2}(y)}$, where $f_{2}(y)=\varphi=\arccos\frac{y(1+\lambda y)}{\lambda-y}$.

From $k=\frac{\pi}{f_{2}(y)}$ and $f_{1}(y)+f(y)=\frac{\pi}{2}$, it follows that \eqref{area} becomes $S=-2\pi\frac{f_{1}(y)}{f_{2}(y)}+2\pi$.

By Lemma \ref{decreasing function}, we get that $\frac{f_{1}(y)}{f_{2}(y)}$ is a decreasing function of $y$. Since $y=g(\frac{\pi}{k})$,
 one can check that $y$ is an increasing function of $k$.
 Consequently, $\frac{f_{1}(y)}{f_{2}(y)}$ is a decreasing function of $k$.

 Because $S=-2\pi\frac{f_{1}(y)}{f_{2}(y)}+2\pi$, the above analysis implies that $S$ is an increasing function of $k$. This completes the proof.
\hfill\qed

\begin{corollary}\label{conj:2}
The area of the regular spherical odd-gon with thickness $\omega\in(0,\frac{\pi}{2})$ is
$2(1-\cos\frac{\omega}{2})\pi$ when the number of vertices tends to infinity.
\end{corollary}

\proof Let $V=v_{1}v_{2}\cdots v_{n}$ be a regular spherical odd-gon.
Since $\lambda=\tan\omega$, we get
\begin{equation*}
\cos\frac{\omega}{2}=\cos\frac{\arctan\lambda}{2}=
 \sqrt{\frac{1+\cos\arctan\lambda}{2}}=\sqrt{\frac{\sqrt{\lambda^{2}+1}+1}{2\sqrt{\lambda^{2}+1}}}.
\end{equation*}
By Theorem \ref{thm3.1}, we get that the area of $V$ is $S=-2\pi\frac{f_{1}(y)}{f_{2}(y)}+2\pi$, where
$y=g(\frac{\pi}{n})=\frac{-(1+\cos\frac{\pi}{n})+\sqrt{(1+\cos\frac{\pi}{n})^{2}+4\lambda^{2}\cos\frac{\pi}{n}}}{2\lambda}$.
When $n$ tends to infinity, $y$ tends to $t=\frac{-1+\sqrt{\lambda^{2}+1}}{\lambda}$.
Then from $\underset{y\to t}{\lim}f_{1}(y)=0$ and $\underset{y\to t}{\lim}f_{2}(y)=0$, by using L$'$Hospital rule, we get
\begin{equation*}
\underset{y\to t}{\lim}
\frac{f_{1}(y)}{f_{2}(y)}
=\underset{y\to t}{\lim}
\frac{f_{1}'(y)}{f_{2}'(y)}
=\underset{y\to t}{\lim}
\frac{\sqrt{1+\lambda^{2}}\sqrt{1+y^{2}}}{-y^{2}+2\lambda y+1}
=\sqrt{\frac{\sqrt{\lambda^{2}+1}+1}{2\sqrt{\lambda^{2}+1}}}.
\end{equation*}
Consequently, $\underset{n\to+\infty }{\lim}S=2(1-\cos\frac{\omega}{2})\pi$.
\hfill\qed

The next theorem shows that the second conjecture mentioned in the introduction is true.

\begin{theorem}\label{thm3.2}
The area of every reduced spherical non-regular $n$-gon
is less than that of the regular spherical $n$-gon of the same thickness.
\end{theorem}

\proof
Let $V=v_{1}\cdots v_{n}$ be a reduced spherical odd-gon.
By Lemma \ref{area of V}, the area of $V$ is
$S=2\sum_{i=1}^n F(\varphi_{i})-(n-2)\pi$, where $F(\varphi_{i})=f(g(\varphi_{i}))=\arcsin\frac{g(\varphi_{i})\sqrt{1+\lambda^{2}}}{\lambda-g(\varphi_{i})}$ and by Fact \ref{fact-varphi}, we have $\varphi_{i}\in(0,\frac{\pi}{2})$.

By Lemma \ref{F(z)}, we obtain that $F(x)$ is a concave function of $x$.
Thus from Jensen's inequality \cite{Jensen}, we have
\begin{equation*}
\frac{F(\varphi_{1})+\cdots+F(\varphi_{n})}{n}\leq F(\frac{\varphi_{1}+\cdots+\varphi_{n}}{n}),
\end{equation*}
the equality holds when $\varphi_{1}=\dots=\varphi_{n}$.
Then the area of $V$ satisfies
\begin{equation*}
S=2n(\frac{F(\varphi_{1})+\cdots+F(\varphi_{n})}{n})-(n-2)\pi
\leq 2n F(\frac{\varphi_{1}+\cdots+\varphi_{n}}{n})-(n-2)\pi.
\end{equation*}

{\bf Case 1.} If $V$ is a regular spherical polygon, then by Claim $2$ in Lemma \ref{sum} we have $\varphi_{1}+\dots+\varphi_{n}=\pi$ and $\varphi_{1}=\dots=\varphi_{n}=\frac{\pi}{n}$. In this case, $S=2n F(\frac{\pi}{n})-(n-2)\pi$.

{\bf Case 2.} If $V$ is a non-regular spherical polygon, then by Claim $1$ in Lemma \ref{sum} we have $\varphi_{1}+\cdots+\varphi_{n}\geq\pi$.

By Lemma \ref{F(z)}, we obtain that $F(x)$ is a decreasing function of $x$.
Since $\varphi_{1}+\dots+\varphi_{n}\geq\pi$, it follows that
$F(\frac{\varphi_{1}+\cdots+\varphi_{n}}{n})\leq F(\frac{\pi}{n})$.
In this case, we have
$S\leq2nF(\frac{\varphi_{1}+\cdots+\varphi_{n}}{n})-(n-2)\pi\leq2n F(\frac{\pi}{n})-(n-2)\pi$.

The above two cases show that the area of $V$ always satisfies $S\leq2n F(\frac{\pi}{n})-(n-2)\pi$, which is exactly the area of the regular spherical $n$-gon.
This completes the proof.
\hfill\qed

By Theorems \ref{thm3.1}, \ref{thm3.2} and Corollary \ref{conj:2} we obtain the following corollary which shows that the first conjecture is true.

\begin{corollary}\label{conj:1}
The area of every reduced spherical polygon $V$ is less than $2(1-\cos\frac{\Delta(V)}{2})\pi$, which is the limit value for the area of the regular spherical odd-gons whose number of vertices tends to infinity.
\end{corollary}


\begin{thebibliography}{n}
\bibitem{Jensen} J. L. W. V. Jensen, Sur les fonctions convexes et les in\'{e}galit\'{e}s
 entre les valeurs moyennes, Acta Math. 1906, 30: 175-193.
\bibitem{Heil} E. Heil, Kleinste konvexe K\"{o}rper gegebener dicke. Preprint No. 453, Fachbereich Mathematik der TH Darmstadt, 1978.
\bibitem{1990} M. Lassak, Reduced convex bodies in the plane. Israel J. Math. 1990, 70: 365-379.
\bibitem{convexbody} M. Lassak, Width of spherical convex bodies. Aequationes Math. 2015, 89: 555-567.
\bibitem{2020} M. Lassak, Diameter, width and thickness of spherical reduced convex bodies with an application to Wulff shapes. Beitr\"{a}ge Alg. Geom. 2020, 61: 369-378.
\bibitem{polygon} M. Lassak, Reduced spherical polygons. Colloq. Math. 2015, 138: 205-216.
\bibitem{2005} M. Lassak, Area of reduced polygons. Publ. Math. 2005, 67: 349-354.
\bibitem{DA} D. A. Murray, Spherical Trigonometry. Longmans Green and CO, London, Bombay and Calcuta 1900.
\bibitem{IP} I. Pak, Lectures on Discrete and Polyhedral Geometry, available at http://www.math.ucla.edu/~pak/geompol8. pdf 2010.
\end{thebibliography}
\end{document}